\documentclass[journal]{IEEEtran}
\IEEEoverridecommandlockouts

\usepackage[latin1]{inputenc}
\usepackage[english]{babel}
\usepackage{amsmath,amsfonts,amssymb, amsthm}
\usepackage{graphicx}
\usepackage{stmaryrd}
\usepackage{epstopdf}
\usepackage{xcolor}
\usepackage{tikz}
\usetikzlibrary{shapes}
\usepackage[]{algorithm2e}
\usepackage{float}
\usepackage{eurosym}

\usepackage{booktabs}

\usepackage{hyperref}

\def\mathscr{\EuScript}

\newcommand{\cF}{\mathcal{F}}

\newcommand{\cU}{\mathcal{U}}

\newcommand{\EE}{\mathbb{E}}

\newcommand{\PP}{\mathbb{P}}

\newcommand{\RR}{\mathbb{R}}

\newcommand{\UU}{\mathbb{U}}

\newcommand{\WW}{\mathbb{W}}
\newcommand{\XX}{\mathbb{X}}

\newcommand{\norm}[1]{\left\|#1\right\|}                    

\newcommand{\opt}{^{\sharp}}                                

\newcommand{\np}[1]{(#1)}                                   
\newcommand{\bp}[1]{\big(#1\big)}                           
\newcommand{\Bp}[1]{\Big(#1\Big)}                           
\newcommand{\bgp}[1]{\bigg(#1\bigg)}                        

\newcommand{\bc}[1]{\big[#1\big]}                           
\newcommand{\Bc}[1]{\Big[#1\Big]}                           

\newcommand{\na}[1]{\{#1\}}                                 
\newcommand{\ba}[1]{\big\{#1\big\}}                         

\newcommand{\argmin}{\mathop{\arg\min}}                     

\makeatletter
\def\va@a{\boldsymbol{\va@arg^{\textstyle\text{\unboldmath$\scriptstyle\va@expo$}}_{\textstyle\text{\unboldmath$\scriptstyle\va@index$}}}}
\def\va#1{\def\va@expo{}\def\va@index{}\def\va@arg{\uppercase{#1}}%
  \@ifnextchar^{\va@h}{\@ifnextchar_\va@u\va@a}}
\def\va@h^#1{\def\va@expo{#1}\@ifnextchar_\va@hu\va@a}
\def\va@u_#1{\def\va@index{#1}\@ifnextchar^\va@uh\va@a}
\def\va@hu_#1{\def\va@index{#1}\va@a}
\def\va@uh^#1{\def\va@expo{#1}\va@a}
\makeatother

\newcommand{\bscal}[2]{\big\langle#1\:,#2\big\rangle}       

\def\eqsepv{\; , \enspace}                                  
\def\eqfinv{\; ,}                                           
\def\eqfinp{\; .}                                           

\newcommand{\finpreuvesymb}{$\Box$}
\newcommand{\finremarksymb}{$\Diamond$}
\newcommand{\finexemplesymb}{$\triangle$}
\newcommand{\finpreuve}{\ \hspace*{\fill}\finpreuvesymb}
\newcommand{\finremark}{\ \hspace*{\fill}\finremarksymb}
\newcommand{\finexemple}{\ \hspace*{\fill}\finexemplesymb}

\makeatletter
\def\endproof{\finpreuve\@endtheorem}
\def\endremark{\finremark\@endtheorem}
\def\endexample{\finexemple\@endtheorem}
\makeatother

\def\x{\boldsymbol{X}}
\def\u{\boldsymbol{U}}
\def\alea{\boldsymbol{W}}
\def\w{\boldsymbol{W}}
\def\fw{\overline{w}}

\def\next{t+1}
\def\post{{\next}}

\def\final{T}

\def\Demandel{\bold{D}^{el}}

\def\Ti{\boldsymbol{\theta}^i}
\def\ti{\theta^i}
\def\Tw{\boldsymbol{\theta}^w}

\def\tank{H}
\def\price{\pi}

\def\opt{^{\sharp}}

\def\price{p}
\def\transfert{\mathbf{F}}
\def\Demandth{\mathbf{D}^{th}}

\def\PV{\boldsymbol{\Phi}^{pv}}

\def\uad{\cU^{ad}}
\def\tw{\theta^w}
\def\to{\theta^e}
\def\pint{\Phi^{int}}
\def\pext{\Phi^{ext}}

\def\tank{h}

\title{Stochastic optimal control of a domestic microgrid \\ equipped with solar panel and battery}

\author{
        Fran\c{c}ois Pacaud\IEEEauthorrefmark{1},
        Pierre Carpentier\IEEEauthorrefmark{2},
        Jean-Philippe Chancelier\IEEEauthorrefmark{3},
        Michel De Lara\IEEEauthorrefmark{3} \\
        \IEEEauthorrefmark{1} Efficacity ---
        \IEEEauthorrefmark{2} ENSTA ---
        \IEEEauthorrefmark{3} CERMICS-ENPC
}

\begin{document}
\maketitle

\begin{abstract}
Microgrids are integrated systems that gather and operate energy production units
to satisfy consumers demands.
This paper details different mathematical methods to design
the Energy Management System (EMS) of domestic microgrids. We
consider different stocks coupled together --- a battery, a domestic hot
water tank --- and decentralized energy production with solar panel.
The main challenge of the EMS is to ensure, at least cost, that supply matches demand for all
time, while considering the inherent uncertainties of such systems.
We benchmark two optimization algorithms to manage the EMS, and compare
them with a heuristic.
The Model Predictive Control (MPC) is a well known algorithm
which models the future uncertainties
with a deterministic forecast. By contrast, Stochastic
Dual Dynamic Programming (SDDP) models the future uncertainties
as probability distributions to compute optimal policies.
We present a fair comparison of these two algorithms to control
microgrid. A comprehensive numerical study shows that
i) optimization algorithms achieve significant gains compared to the heuristic,
ii) SDDP outperforms MPC by a few percents, with a reasonable computational overhead.
\end{abstract}



\section{Introduction}

\subsection{Context}

A microgrid is a local energy network that produces part of its energy
and controls its own demand. Such systems are complex to control, because
of the different stocks and interconnections.
Furthermore, at local scale, electrical demands and weather conditions
(heat demand and renewable energy production) are highly variable and
hard to predict; their stochastic nature adds uncertainty to the system.

\begin{figure}[!ht]
    \begin{center}
      \input{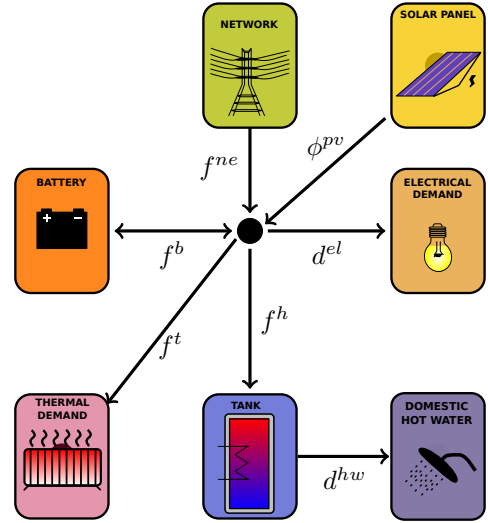}
    \end{center}
    \caption{Electrical microgrid }
    \label{fig:twostocksschema}
\end{figure}

We consider here a domestic microgrid (see Figure~\ref{fig:twostocksschema}),
equipped with a battery, an electrical hot water tank and a solar panel.
We use the battery to store energy when prices are low or when the
production of the solar panel is above the electrical demand.
The microgrid is connected to an external grid to import electricity when
needed.
Furthermore, we model the building's envelope to take advantage of the thermal inertia
of the building.
Hence, the system has four stocks to store energy: a battery, a hot water tank,
and two passive stocks being the building's walls and inner rooms.
Two kind of uncertainties affect the system. First, the electrical
and domestic hot water demands are not known in advance.
Second, the production of the
solar panel is heavily perturbed because of the varying nebulosity affecting
their production.

We aim to compare two classes of algorithms to tackle the uncertainty
in microgrid Energy Management Systems (EMS). The renowned Model Predictive
Control (MPC) algorithm views the future uncertainties with a deterministic
forecast. Then, MPC relies on deterministic optimization algorithms to compute
optimal decisions. The contender --- Stochastic
Dual Dynamic Programming (SDDP) --- is an algorithm based on the
Dynamic Programming principle. Such algorithm computes offline a set of value functions
by backward induction; optimal decisions are computed online as time goes on,
using the value functions.
We present a balanced comparison of these two algorithms, and highlight the
advantages and drawbacks of both methods.

\subsection{Litterature}

\subsubsection{Optimization and energy management systems}
Energy Management Systems (EMS) are integrated automated tools used to
monitor and control energy systems.
The design of EMS for buildings has raised interest in
recent years.
In~\cite{olivares2014trends}, the authors give an overview concerning
the application of optimization methods in designing EMS.

The well-known Model Predictive Control (MPC)~\cite{garcia1989model} has been
widely used to control EMS.
We refer notably to \cite{oldewurtel2011stochastic}, \cite{malisani2012thesis},
\cite{lamoudi2012distributed} for applications of
MPC in buildings.
Different solutions are investigated to tackle uncertainties,
such as Stochastic MPC~\cite{oldewurtel2011stochastic}
or robust optimization~\cite{paridari2016robust}.

\subsubsection{Stochastic Optimization}
as we said, at local scale, electrical demand and production are highly
variable, especially as microgrids are expected to absorb renewable energies.
This leads to pay attention to stochastic optimization approaches.
Apart from microgrid management, stochastic optimization has found some applications
in energy systems (see~\cite{delara2014cfe} for an overview).
Historically, stochastic optimization has been
widely applied to hydrovalleys management \cite{pereira1991multi}.
Other applications have arisen recently, such as integration of wind energy
and storage \cite{haessig2014dimensionnement} or insulated microgrids management
\cite{heymann2016stochastic}.

Stochastic Dynamic Programming (SDP)~\cite{bertsekas1995dynamic} is a general method
to solve stochastic optimal control problems.
In energy applications, a variant of SDP, Stochastic Dual Dynamic Programming (SDDP),
has proved its adequacy for large scale applications.
SDDP was first described in the seminal paper~\cite{pereira1991multi}. We
refer to~\cite{shapiro2011analysis} for a generic
description of the algorithm and its application to the management of hydrovalleys.
A proof of convergence in the linear case is given in
\cite{philpott2008convergence}, and in the convex case
in~\cite{girardeau2014convergence}.

With the growing adoption of stochastic optimization methods, new researches
aim to compare algorithms such as SDP and SDDP with
MPC. We refer to the recent paper~\cite{riseth2017comparison}.

\subsection{Structure of the paper}

We detail a modelling of the microgrid in Sect.~\ref{sec:problem}, then
formulate an optimization problem in Sect.~\ref{sec:optimization}.
We outline the different optimization algorithms in Sect.~\ref{sec:algo}.
Finally, we provide in Sect.~\ref{sec:numeric} numerical results concerning the management
of the microgrid.

\section{Energy system model}
\label{sec:problem}

In this section, we depict the physical equations of the energy system model described
in Figure~\ref{fig:twostocksschema}. These equations write
naturally in continuous time~$t$.
We model the battery and the hot water tank with stock dynamics, and
describe the dynamics of the building's temperatures with an electrical analogy.
Such physical model fulfills two purposes:
it will be used to assess different control policies;
it will be the basis of a discrete time model used to design optimal control policies.

\subsection{Load balance}

Based on Figure~\ref{fig:twostocksschema}, the \emph{load
balance} equation of the microgrid writes, at each time~$t$:
\begin{equation}
  \label{eq:loadbalance}
  \phi^{pv}(t) + f^{ne}(t) = f^b(t) + f^t(t) + f^h(t) + d^{el}(t) \eqfinp
\end{equation}
We now comment the different terms.
In the left hand side of Equation~\eqref{eq:loadbalance}, the load produced consists of
\begin{itemize}
  \item the production of the solar panel~$\phi^{pv}(t)$,
  \item the importation from the network~$f^{ne}(t)$, supposed nonnegative (we do not
    export electricity to the network).
\end{itemize}
In the right hand side of Equation~\eqref{eq:loadbalance}, the electrical demand
is the sum of
\begin{itemize}
  \item the power exchanged with the battery~$f^b(t)$,
  \item the power injected in the electrical heater~$f^t(t)$,
  \item the power injected in the electrical hot water tank~$f^h(t)$,
  \item the inflexible demands (lightning, cooking...), aggregated in a single
    demand~$d^{el}(t)$.
\end{itemize}

\subsection{Energy storage}
\label{sec:twostockselec}

We consider a lithium-ion battery, whose state of charge at time~$t$ is
denoted by~$b(t)$. The state of charge is bounded:
\begin{equation}
  \label{eq:batterybound}
  \underline b \leq b(t) \leq \overline b \eqfinp
\end{equation}
Usually, we set~$\underline b = 30\% \times \overline b$ so as to avoid
empty state of charge, which proves to be stressful for the battery.
The battery dynamics is given by the differential equation
\begin{equation}
  \label{eq:batteryequation}
  \dfrac{db}{dt} = \rho_c (f^b(t))^+ - \dfrac{1}{\rho_d} (f^b(t))^- \eqfinv
\end{equation}
with~$\rho_c$ and~$\rho_d$ being the charge and discharge efficiency
and~$f^b(t)$ denoting the power exchange with the battery. We use the convention
$f^{+} = \max\np{0, f}$ and $f^{-} = \max\np{0, -f}$.

As we cannot withdraw an infinite power from the battery at time~$t$,
we bound the power exchanged with the battery:
\begin{equation}
  - \underline f^b \leq f^b(t) \leq \overline f^b \eqfinp
\end{equation}

\subsection{Electrical hot water tank}

We use a simple linear model for the electrical hot water tank dynamics.
At time~$t$, we denote by~$T^h(t)$ the temperature inside the hot water tank.
We suppose that this temperature is homogeneous, that is, that no
stratification occurs inside the tank.

At time~$t$, we define the energy~$\tank(t)$ stored inside the tank
as the difference between the tank's temperature
$T^{\tank}(t)$ and a reference temperature~$T^{ref}$
\begin{equation}
  \tank(t) = \rho V_{\tank} c_p \bp{T^{\tank}(t) - T^{ref}} \eqfinv
\end{equation}
where~$V_{\tank}$ is the tank's volume,~$c_p$ the calorific capacity of water and
$\rho$ the density of water.
The energy~$\tank(t)$ is bounded:
\begin{equation}
  \label{eq:tankbounds}
  0 \leq \tank(t) \leq \overline h  \eqfinp
\end{equation}
The enthalpy balance equation writes
\begin{equation}
    \label{eq:tankequation}
    \dfrac{dh}{dt} = \beta_\tank f^h(t) - d^{hw}(t) \eqfinv
\end{equation}
where
\begin{itemize}
  \item $f^{h}(t)$ is the electrical power used to heat the tank, satisfying
    \begin{equation}
    0 \leq f^h(t) \leq \overline f^h \eqfinv
    \end{equation}
  \item $d^{hw}(t)$ is the domestic hot water demand,
  \item $\beta_\tank$ is a conversion yield.
\end{itemize}
A more accurate representation would model the stratification inside the hot water
tank. However, this would greatly increase the number of states in the system,
rendering the numerical resolution more cumbersome. We refer to \cite{schutz2015comparison}
and \cite{beeker2016discrete} for discussions about the impact of the tank's
modeling on the performance of the control algorithms.

\subsection{Thermal envelope}
We model the evolution of the temperatures inside the building with
an electrical analogy: we view temperatures as voltages, walls as capacitors,
and thermal flows as currents. A model with 6 resistances and 2 capacitors
(R6C2) proves to be accurate to describe small buildings \cite{berthou2013thesis}.
The model takes into account two temperatures:
\begin{itemize}
  \item the wall's temperature~$\tw(t)$,
  \item the inner temperature~$\ti(t)$.
\end{itemize}
Their evolution is governed by the two following differential
equations
\label{sec:twostockstherm}
\begin{subequations}
  \label{eq:r6c2continuous}
  \begin{multline}
  c_m \dfrac{d \tw}{dt} = \underbrace{\dfrac{\ti(t) - \tw(t)}{R_i + R_s}}_{\substack{\text{Exchange} \\ \text{  Indoor/Wall}}}
                              + \underbrace{\dfrac{\to(t) - \tw(t)}{R_m + R_e}}_{\substack{\text{Exchange} \\ \text{  Outdoor/Wall}}}
                              + \underbrace{\gamma f^{t}(t)}_{\text{Heater}}\\
                               + \underbrace{\dfrac{R_i}{R_i + R_s} \pint(t)}_{\substack{\text{Radiation} \\ \text{through windows}}}
                              + \underbrace{\dfrac{R_e}{R_e + R_m} \pext(t)}_{\substack{\text{Radiation} \\ \text{through wall}}} \eqfinv
  \end{multline}
\begin{multline}
  c_i \dfrac{d \ti}{dt} = \underbrace{\dfrac{\tw(t) - \ti(t)}{R_i + R_s}}_{\substack{\text{Exchange} \\ \text{  Indoor/Wall}}}
                              + \underbrace{\dfrac{\to(t) - \ti(t)}{R_v}}_{\text{Ventilation}}
                              + \underbrace{\dfrac{\to(t) - \ti(t)}{R_f}}_{\text{Windows}}\\
                              + \underbrace{(1 - \gamma) f^{t}(t)}_{\text{Heater}}
                              + \underbrace{\dfrac{R_s}{R_i + R_s} \pint(t)}_{\substack{\text{Radiation} \\ \text{through windows}}} \eqfinv
\end{multline}
\end{subequations}
where we denote
\begin{itemize}
  \item the power injected in the heater by~$f^t(t)$,
  \item the external temperature by~$\to(t)$,
  \item the radiation on the wall by~$\pext(t)$,
  \item the radiation through the windows by~$\pint(t)$.
\end{itemize}
The time-varying quantities~$\to(t)$, $\pint(t)$ and~$\pext(t)$ are exogeneous.
We denote by~$R_i, R_s, R_m, R_e, R_v, R_f$ the different resistances of the
R6C2 model,
and by~$c_i, c_m$ the capacities of the inner rooms
and the walls. We denote by $\gamma$ the proportion of heating dissipated in the wall
through conduction, and by $(1- \gamma)$ the proportion of heating dissipated
in the inner room through convection.

\subsection{Continuous time state equation}

We denote by~$x = (b, h, \tw, \ti)$ the state,
$u = (f^b, f^t, f^h)$ the control, and~$w = (d^{el}, d^{hw},\phi^{pv})$ the
uncertainties. The continuous state equation writes
\begin{equation}
  \label{eq:continuousSE}
  \dot{x} = F(t, x, u, w) \eqfinv
\end{equation}
where the function~$F$ is defined by
Equations~\eqref{eq:batteryequation}-\eqref{eq:tankequation}-\eqref{eq:r6c2continuous}.

\section{Optimization problem statement}
\label{sec:optimization}

Now that we have described the physical model, we turn to the formulation of
a decision problem. We aim to compute optimal decisions that
minimize a daily operational cost,
by solving a stochastic optimization problem.

\subsection{Decisions are taken at discrete times}

The EMS takes decisions every 15 minutes to control the system.
Thus, we have to provide decisions in discrete time.

We set~$\Delta = 15 \text{mn}$, and we consider an horizon~$\final_0 = 24\text{h}$.
We adopt the following convention
for discrete processes:  for~$t \in \{0, 1, \cdots, \final = \frac{\final_0}{\Delta}\}$, we set
$x_t = x(t \Delta)$. That is,~$x_t$ denotes the value of the variable~$x$ at the
beginning of the interval~$[t \Delta, (t+1) \Delta[$. Otherwise stated, we will
denote by $[t, t+1[$ the continuous time interval $[t \Delta, (t+1) \Delta[$.

\subsection{Modeling uncertainties as random variables}

Because of their unpredictable nature, we cannot anticipate
the realizations of the electrical and the thermal demands.
A similar
reasoning applies to the production of the solar panel.
We choose to model these quantities as random variables (over a sample space~$\Omega$).
We adopt the following convention: a random variable will be denoted by an uppercase bold letter~$\va Z$
and its realization will be denoted in lowercase~$z = \va Z(\omega)$.
For each~$t = 1, \ldots, \final$, we define the uncertainty vector
\begin{equation}
  \w_t = (\Demandel_t, \Demandth_t, \PV_t) \eqfinv
\end{equation}
modeled as a random variable. The uncertainty
$\w_t$ takes value in the set~$\WW_t = \RR^3$.

\subsection{Modeling controls as random variables}

As decisions depend on the previous uncertainties, the control is a random
variable.
We recall that, at each discrete time~$t$, the EMS takes three decisions:
\begin{itemize}
  \item how much energy to charge/discharge the battery~$\va f_t^b$,
  \item how much energy to store in the electrical hot
    water tank~$\va f_t^h$,
  \item how much energy to inject in the electrical heater
      $\va f_t^t$.
\end{itemize}
We write the decision vector (random variable)
\begin{equation}
  \u_t = (\va f_t^b, \va f_t^h, \va f_t^t) \eqfinv
\end{equation}
taking values in~$\UU_t = \RR^3$.

Then, between two discrete time indexes~$t$ and~$t+1$, the EMS imports an energy~$\va
f^{ne}_\post$ from the external network. The EMS must fulfill the load balance
equation~\eqref{eq:loadbalance}
whatever the demand~$\Demandel_\post$ and the production of the solar panel
$\PV_\post$, unknown at time~$t$. Hence~$\va f^{ne}_\post$ is a recourse
decision taken at time~$\post$.
The load balance equation~\eqref{eq:loadbalance} now writes as
\begin{equation}
  \label{eq:loadbalancerandom}
  \va f_\post^{ne}= \va f_t^b + \va f^t_t + \va f^h_t + \Demandel_\post-\PV_\post
  \quad \PP-\text{a.s.} \eqfinv
\end{equation}
where~$\PP-\text{a.s.}$ indicates that the constraint is fulfilled in the
almost sure sense.
Later, we will aggregate the solar panel production~$\PV_\post$ with
the demands~$\Demandel_\post$ in Equation~\eqref{eq:loadbalancerandom},
as these two quantities appear only by their sum.

The need of a recourse variable is a consequence of stochasticity in the
supply-demand equation.
The choice of the recourse variable depends on the modeling.
Here, we choose the recourse $\va f^{ne}_\post$ to be
provided by the external network, that is, the external network mitigates
the uncertainties in the system.

\subsection{States and dynamics\label{sec:state}}
The state becomes also a random variable
\begin{equation}
  \x_t = (\va b_t, \va h_t, \Tw_t, \Ti_t) \eqfinp
\end{equation}
It gathers the stocks in the battery~$\va b_t$ and
in the electrical hot water tank~$\va h_t$, plus the two temperatures of
the thermal envelope~$(\Ti_t, \Tw_t)$. Thus, the state vector $\x_t$
takes values in~$\XX_t = \RR^4$.

The discrete dynamics writes
\begin{equation} \label{eq:dynamics}
  x_\post = f_t(x_t, u_t, w_\post)\eqfinv
\end{equation}
where~$f_t$ corresponds to the discretization of the continuous
dynamics~\eqref{eq:continuousSE} using a forward Euler scheme, that is,
$x_\post = x_t + \Delta  \times  F(t, x_t, u_t, w_\post)$.
By doing so, we suppose that the control~$u_t$ and the uncertainty~$w_\post$
are constant over the interval~$[t, t+\Delta[$.

The dynamics~\eqref{eq:dynamics} rewrites as an almost-sure constraint:
\begin{equation}
    \x_{t+1} = f_t \bp{\x_t, \u_t, \alea_{t+1}} \qquad \PP-\text{a.s.} \eqfinp
\end{equation}
We suppose that we start from a given position~$x_0$, thus adding a new
initial constraint:~$\x_0 = x_0$.

\subsection{Non-anticipativity constraints}

The future realizations of uncertainties are unpredictable.
Thus, decisions are functions of previous history only, that is,
the information collected between time~$0$ and time~$t$. Such a constraint
is encoded as an algebraic constraint, using the tools of Probability
theory~\cite{kallenberg2002}.
The so-called non-anticipativity constraint writes
\begin{equation}
  \label{eq:measurability}
  \sigma(\u_t) \subset \cF_t \eqfinv
\end{equation}
where~$\sigma(\u_t)$ is the~$\sigma$-algebra generated by~$\u_t$ and
$\cF_t = \sigma(\w_1, \cdots, \w_t)$ the~$\sigma$-algebra
associated to the previous history~$(\w_1, \ldots, \w_t)$.
If Constraint~\eqref{eq:measurability} holds true, the Doob lemma~\cite{kallenberg2002}
ensures that there exists a function~$\pi_t$ such that
\begin{equation}
  \u_t = \pi_t(\x_0, \w_1, \ldots, \w_t) \eqfinp
\end{equation}
This is how we turn an (abstract) algebraic constraint into a more practical
functional constraint. The function~$\pi_t$ is an example of policy.

\subsection{Bounds constraints}

By Equations~\eqref{eq:batterybound} and \eqref{eq:tankbounds}, the stocks in the battery~$\va b_t$ and in the tank~$\va h_t$ are bounded.
At time~$t$, the control~$\va f^b_t$ must ensure that the next state~$\va b_\post$
is admissible, that is, $\underline b \leq \va b_\post \leq \overline b$
by Equation~\eqref{eq:batterybound}, which rewrites,
\begin{equation}
  \underline b \leq
  \va b_t + \Delta \bc{\rho_c (\va f_t^b)^+ + \dfrac{1}{\rho_d} (\va f_t^b)^-}
  \leq \overline b \eqfinp
\end{equation}
Thus, the constraints on~$\va f_t^b$ depends on the stock~$\va b_t$.
The same reasoning applies for the tank power~$\va f_t^h$.
Furthermore, we set bound constraints on controls, that is,
\begin{equation}
  -\overline f^b \leq \va f_t^b \leq \overline f^b \eqsepv
  0 \leq \va f_t^h \leq \overline f_t^h \eqsepv
  0 \leq \va f_t^t \leq \overline f_t^t \eqfinp
\end{equation}
Finally the load-balance equation~\eqref{eq:loadbalancerandom} also acts as a constraint
on the controls.
We gather all these constraints in an admissible set on control~$\u_t$
depending on the current state~$\x_t$:
\begin{equation}
  \u_t \in \uad_t(\x_t) \qquad \PP-\text{a.s.} \eqfinp
\end{equation}

\subsection{Objective}
\label{sec:objective}
At time $t$, the operational cost~$L_t: \XX_t \times \UU_t \times \WW_\post \rightarrow \RR$
aggregates two different costs:
\begin{equation}
  \label{eq:operational_cost}
  L_t(x_t, u_t, w_\post) = \price^e_t  \times f^{ne}_\post +
\price^d_t\times  \max\np{0, \overline{\theta^i_t} - \ti_t} \eqfinp
\end{equation}
First, we pay a price~$\price_t^e$ to import electricity from the network
between time~$t$
and~$t+1$. Hence, electricity cost is equal
to~$\price^e_t \times \va f^{ne}_\post  $.
Second, if the indoor temperature is below a given threshold, we penalize the
induced discomfort with a cost
  $\price^d_t \times \max\np{0, \overline{\theta^i_t} - \Ti_t}$,
where~$\price_t^d$ is a virtual price of discomfort.
The cost~$L_t$ is a convex piecewise linear function, which will prove important
for the SDDP algorithm.

We add a final cost~$K : \XX_{\final} \rightarrow \RR$
to ensure that stocks are non empty at final
time~$\final$
\begin{equation}
  \label{eq:finalpenal}
  K(x_\final) = \kappa \times \max \np{0, x_0 - x_\final} \eqfinv
\end{equation}
where~$\kappa$ is a positive penalization coefficient calibrated by trials
and errors.

As decisions $\u_t$ and states $\x_t$ are random, the costs $L_t(\x_t, \u_t, \w_\post)$
become also random variables. We choose to minimize the expected value
of the daily operational cost, yielding the criterion
\begin{equation}
  \label{eq:objective}
  \EE \Bc{\sum_{t=0}^{\final - 1} L_t(\x_t, \u_t, \w_\post)
  + K(\x_\final) } \eqfinv
\end{equation}
yielding an expected value of a convex piecewise linear cost.

\subsection{Stochastic optimal control formulation}
Finally, the EMS problem translates to a generic Stochastic Optimal Control (SOC)
problem
\begin{subequations}
  \label{eq:stochproblem}
\begin{align}
  \underset{\x, \u}{\min} ~ & \EE \; \Bc{\sum_{t=0}^{\final-1}
                              L_t(\x_t, \u_t, \alea_{t+1}) + K(\x_\final)} \eqfinv \\
                            & \x_0 = x_0  \eqfinv\\
                            & \x_{t+1} = f_t \bp{\x_t, \u_t, \alea_{\next}}\quad \PP-\text{a.s.} \eqfinv \\
                            & \u_t \in \uad_t(\x_t)  \quad \PP-\text{a.s.} \eqfinv \\
                            & \sigma(\u_t) \subset \cF_t \eqfinp
\end{align}
\end{subequations}
Problem \eqref{eq:stochproblem} states that we want to minimize the expected
value of the costs while satisfying the dynamics, the control bounds and
the non-anticipativity constraints.

\section{Resolution methods}
\label{sec:algo}

The exact resolution of Problem~\eqref{eq:stochproblem} is out of reach in
general. We propose two different algorithms that provide policies
$\pi_t: \XX_0 \times \WW_1 \times \cdots \times \WW_t \rightarrow \UU_t$
that map available information~$x_0, w_1, \ldots, w_t$ at time~$t$ to a decision~$u_t$.

\subsection{Model Predictive Control (MPC)}
\label{MPC}

MPC is a classical algorithm commonly used to handle uncertainties in energy
systems.
At time~$t$, it considers a deterministic forecast~$(\fw_\post, \ldots,
\fw_\final)$ of the future uncertainties~$(\w_\post, \ldots, \w_\final)$ and
solves the deterministic problem
\begin{subequations}
  \label{eq:mpcproblem}
\begin{align}
   \underset{(u_t, \cdots, u_{\final-1})}{\min} ~ & \sum_{j=t}^{\final-1}\Bc{
                              L_j(x_j, u_j, \fw_{j+1})} + K(x_\final)\eqfinv \\
                            & x_{j+1} = f_j \bp{x_j, u_j, \fw_{j+1}} \eqfinv\\
                            & u_j \in \uad_j(x_j) \eqfinp
\end{align}
\end{subequations}
At time $t$, we solve Problem~\eqref{eq:mpcproblem}, retrieve the
optimal decisions~$(u_t^{\opt}, \ldots, u_{\final-1}^{\opt})$ and only keep the
first decision~$u_t^{\opt}$ to control the system between time~$t$ and~$t+1$.
Then, we restart the procedure at time $t+1$.

As Problem~\eqref{eq:mpcproblem} is linear and the number of time steps
remains not too large, we are able to solve it exactly for every~$t$.

\subsection{Stochastic Dual Dynamic Programming (SDDP)}
\subsubsection{Dynamic Programming and Bellman principle}
the Dynamic Programming method~\cite{Bellman57} looks for solutions of
Problem~\eqref{eq:stochproblem}
as state-feedbacks $\pi_t: \XX_t \rightarrow \UU_t$. Dynamic Programming
computes a serie of value functions
backward in time by setting $V_\final(x_\final) = K(x_\final)$ and solving
the recursive equations
\begin{multline}
\label{eq:bellmanequations}
  V_t(x_t) = \min_{u \in \uad_t(x_t)}
  \int_{\WW_\post} \Bc{ L_t(x_t, u, w_\post) + \\
  V_\post\bp{f_t(x_t, u, w_\post)}}
  \mu^{of}_\post(dw_\post) \eqfinv
\end{multline}
where~$\mu^{of}_\post$ is a (offline) probability distribution on~$\WW_\post$.

Once these functions are computed, we compute a decision at time $t$
as a state-feedback:
\begin{multline}
  \label{eq:bellmanpolicy}
  \pi_t(x_t) \in \argmin_{u \in \uad_t(x_t)}
  \int_{\WW_\post} \Bc{ L_t(x_t, u, w_\post) + \\
  V_\post\bp{f_t(x_t, u, w_\post)}}
  \mu^{on}_\post(dw_\post) \eqfinv
\end{multline}
where $\mu_\post^{on}$ is an online probability distribution on~$\WW_\post$.
This method proves to be optimal when the uncertainties~$\w_1, \ldots, \w_\final$
are stagewise independent and when $\mu^{on}_t = \mu^{of}_t$ is the probability distribution
of~$\w_t$
in~\eqref{eq:bellmanequations}.

\subsubsection{Description of Stochastic Dual Dynamic Programming}
\label{sec:sddp}
Dynamic Programming suffers from the well-known \emph{curse
of dimensionality}~\cite{bertsekas1995dynamic}: its numerical resolution fails
for state dimension typically greater than~4 when value functions are computed
on a discretized grid.
As the state~$\x_t$ in~\S\ref{sec:state} has 4 dimensions,
SDP would be too slow to solve numerically Problem~\eqref{eq:stochproblem}.
The Stochastic Dual Dynamic Programming (SDDP) can bypass
the curse of dimensionality by approximating value functions by polyhedral
functions. It is optimal for solving Problem~\eqref{eq:stochproblem} when
uncertainties are stagewise independent, costs~$L_t$ and~$K$ are convex and
dynamics~$f_t$ are linear~\cite{girardeau2014convergence}.

SDDP provides an outer approximation~$\underline{V}_t^k$ of the value function
$V_t$ in \eqref{eq:bellmanequations}
with a set of supporting hyperplanes
$\na{(\lambda^j_t, \beta^j_t)}_{j = 1, \cdots, k}$ by
\begin{subequations}
  \begin{align}
    \underline{V}_t(x_t) &= \min_{\theta_t \in \RR} \theta_t \eqfinv\\
                         & \bscal{\lambda_t^j}{x_t} + \beta_t^j \leq \theta_t
    \eqsepv \forall j = 1, \cdots, k \eqfinp
  \end{align}
\end{subequations}
Each iteration~$k$ of SDDP encompasses two passes.
\begin{itemize}
  \item During the \emph{forward pass}, we draw a scenario~$x_0, \ldots,
    w_\final^k$
    of uncertainties,
    and compute a state trajectory~$\ba{x_t^k}_{t = 0 \cdots \final}$
    along this scenario.
    Starting from position~$x_0$, we compute~$x_\post^k$ in an iterative
    fashion: i) we compute the optimal control at time~$t$ using the available
    $\underline V_\post^k$ function
    \begin{multline}
      \label{eq:sddpforward}
      u_t^k \in \argmin_{u \in \uad_t(x_t)} \int_{\WW_\post} \bc{
      L_t(x_t^k, u, w_\post) + \\ \underline V_\post^k\bp{f_t(x_t^k, u, w_\post)}}
      \mu_\post^{of}(dw_\post)
      \eqfinv
    \end{multline}
    and ii), we set~$x_\post^k = f_t(x_t^k, u_t^k, w_\post^k)$ where $f_t$ is
    given by~\eqref{eq:dynamics}.
  \item During the \emph{backward pass}, we update the approximated value functions
    $\ba{\underline V_t^k}_{t=0, \cdots, \final}$ backward in time along the
    trajectory~$\ba{x_t^k}_{t=0, \cdots, \final}$. At time $t$,
    we solve the problem
    \begin{multline}
      \label{eq:sddpbackward}
      \theta_t^{k+1} =  \min_{u \in \uad_t(x_t)} \int_{\WW_\post}\bc{
      L_t(x_t^k, u, w_\post) + \\ \underline V_\post^{k+1} \bp{f_t(x_t^k, u, w_\post)}}
      \mu_\post^{of}(dw_\post)
      \eqfinv
    \end{multline}
    and we obtain a new cut~$(\lambda_t^{k+1}, \beta_t^{k+1})$ where~$\lambda_t^{k+1}$
    is a subgradient of optimal cost~\eqref{eq:sddpbackward} evaluated at point $x_t = x_t^k$ and
    $\beta_t^{k+1} = \theta_t^{k+1} - \bscal{\lambda_t^{k+1}}{x_t^{k}}$.
    This new cut allows to update the function~$\underline V_t^{k+1}$:
    $V_t^{k+1} = \max\na{\underline V_t^k, \bscal{\lambda_t^{k+1}}{.} + \beta_t^{k+1}}$.
\end{itemize}
Otherwise stated, SDDP only explores
the state space around ``interesting'' trajectories (those computed during
the forward passes) and refines the value functions only in the corresponding
space regions (backward passes).

\subsubsection{Obtaining online controls with SDDP}
in order to compute implementable decisions, we use the following procedure.
\begin{itemize}
  \item Approximated value functions~$\ba{\underline V_t}$ are computed
  with the SDDP algorithm (see~\S\ref{sec:sddp}). These computations are
  done offline.
  \item The approximated value functions~$\ba{\underline V_t}$ are then used
  to compute online a decision at any time~$t$ for any state~$x_t$.
\end{itemize}
More precisely, we compute the SDDP policy~$\pi^{sddp}_t$ by
\begin{multline}
  \label{eq:sddppolicy}
  \pi^{sddp}_t(x_t) \in \argmin_{u \in \uad_t(x_t)} \int_{\WW_\post}\bc{
    L_t(x_t, u, w_\post) + \\ \underline V_\post\bp{f_t(x_t, u, w_\post)}}
    \mu_\post^{on}(dw_\post)
    \eqfinv
\end{multline}
which corresponds to replacing the value function~$V_\post$ in
Equation~\eqref{eq:bellmanpolicy} with its approximation $\underline V_\post$.
The decision $\pi_t^{sddp}(x_t)$ is used to control the system between time~$t$
and~$t+1$. Then, we resolve Problem~\eqref{eq:sddppolicy} at time $t+1$.

To solve numerically problems~\eqref{eq:sddpforward}-\eqref{eq:sddpbackward}-\eqref{eq:sddppolicy}
at time~$t$,
we will consider distributions with finite support $w_t^1, \ldots, w_t^S$.
The offline distribution~$\mu_\post^{of}$ now writes: $\mu_\post^{of} = \sum_{s=1}^S p_s \delta_{w_\post^s}$
where $\delta_{w_\post^s}$ is the Dirac measure at~$w_\post^s$
and $(p_1, \ldots, p_S)$ are probability weights. The same reasoning applies
to the online distribution~$\mu_\post^{on}$. For instance,
Problem~\eqref{eq:sddppolicy} reformulates as
\begin{multline}
  \label{eq:sddppolicydiscrete}
  \pi^{sddp}_t(x_t) \in \argmin_{u \in \uad_t(x_t)} \sum_{s=1}^S p_s \bc{
    L_t(x_t, u, w_\post^s) + \\ \underline V_\post\bp{f_t(x_t, u, w_\post^s)}}
    \eqfinp
\end{multline}

\section{Numerical results}
\label{sec:numeric}

\subsection{Case study}
\subsubsection{Settings}
we aim to solve the stochastic optimization problem~\eqref{eq:stochproblem}
over one day, with 96 time steps. The battery's size is~$3~\text{kWh}$, and the
hot water tank has a capacity of~$120~\text{l}$.
We suppose that the house has a surface~$A_p = 20~\text{m}^2$ of solar panel
at disposal, oriented south, and with a yield of
15\%. We penalize the recourse variable
$\transfert^{ne}_{t+1}$ in~\eqref{eq:operational_cost} with on-peak and off-peak tariff, corresponding to
\'Electricit\'{e} de France's (EDF) individual tariffs.
The building's thermal envelope corresponds
to the French RT2012 specifications~\cite{rt2012}.
Meteorological data comes from Meteo France measurements corresponding
to the year 2015.

\subsubsection{Demands scenarios}
\begin{figure}[!ht]
  \begin{tabular}{cc}
    \includegraphics[width=4.5cm]{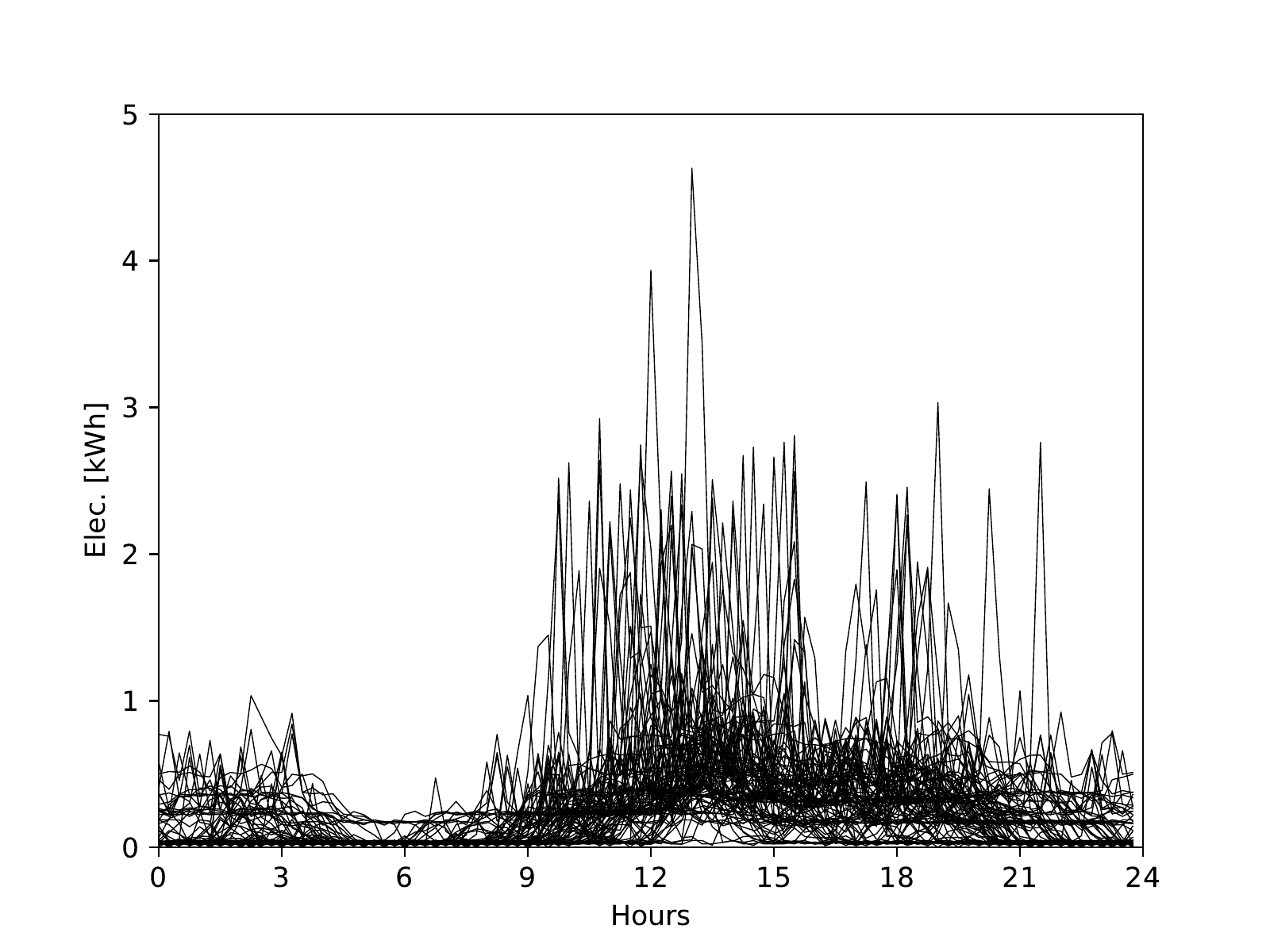} & \hspace{-0.6cm}
    \includegraphics[width=4.5cm]{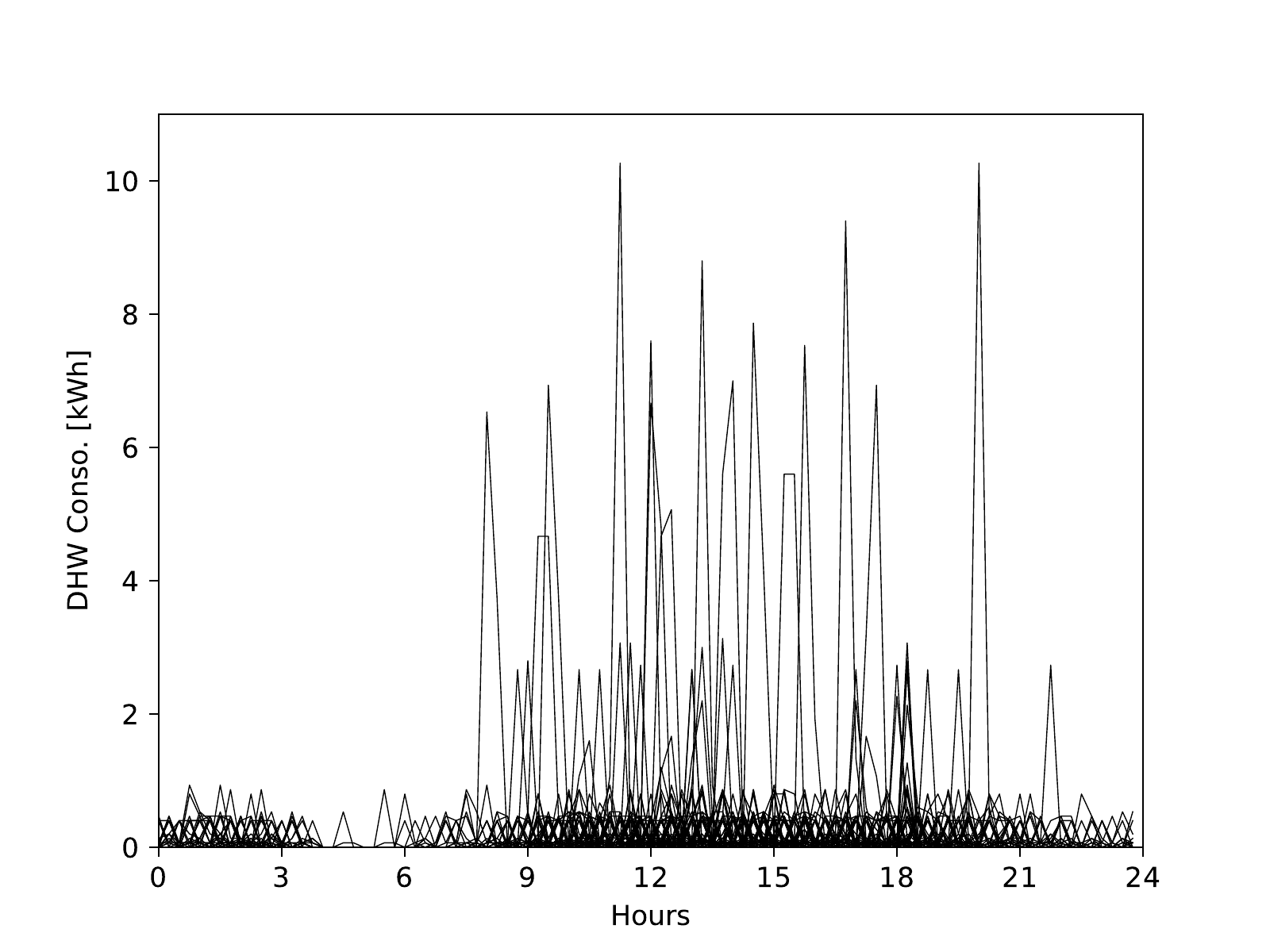}
  \end{tabular}
  \caption{Electrical (left) and domestic hot water (right) demand scenarios.}
  \label{fig:scenstrobe}
\end{figure}

we have scenarios of electrical and domestic hot water demands at 15 minutes
time steps, obtained with StRoBe~\cite{baetens2016modelling}.
Figure~\ref{fig:scenstrobe} displays 100 scenarios of electrical and hot
water demands.
We observe that the shape of these scenarios is consistent: demands are almost
null during night, and there are peaks around midday and 8~pm.
Peaks in hot water demands corresponds to showers.
We aggregate the production of the solar panel~$\PV$ and the electrical demands
$\Demandel$ in a single variable $\Demandel$ to consider
only two uncertainties~$(\Demandel_t, \va d^{hw}_t)$.

\subsubsection{Out of sample assessment of strategies}
to obtain a fair comparison between SDDP and MPC, we use an out-of-sample
validation.
We generate 2,000 scenarios of electrical and hot water demands,
and we split these scenarios in two distinct parts: the first~$N_{opt}$
scenarios are called \emph{optimization scenarios}, and the
remaining~$N_{sim}$ scenarios are called \emph{assessment scenarios}.
We made the choise $N_{opt} = N_{sim} = 1,000$.

First, during the offline phase, we use the
optimization scenarios to build models for the uncertainties,
under the mathematical form required by each algorithm
(see Sect.~\ref{sec:algo}).
Second, during the online phase, we use the assessment scenarios to compare
the strategies produced by these algorithms.
At time~$t$ during the assessment, the algorithms cannot use the future
values of the assessment scenarios, but can take advantage
of the observed values up to~$t$ to update their statistical models of future uncertainties.

\subsection{Numerical implementation}

\subsubsection{Implementing the algorithms}
we implement MPC and SDDP in Julia 0.6, using JuMP~\cite{dunning2017jump}
as a modeler, \verb+StochDynamicProgramming.jl+ as a SDDP solver, and
Gurobi 7.02~\cite{gurobi2015gurobi} as a LP solver.
All computations run on a Core i7 2.5 GHz processor, with 16Go RAM.

\subsubsection{MPC procedure}
Electrical and thermal demands are naturally correlated in
time~\cite{widen2010high}.
To take into account such a dependence across the different
time-steps, we chose to model the process~$\w_1,
\ldots, \w_\final$ with an auto-regressive (AR) process.

\paragraph{Building offline an AR model for MPC}
we fit an AR(1) model upon the optimization scenarios
(we do not consider higher order lag for the sake of simplicity).
For $i \in \{el, hw\}$, the AR model writes
\begin{subequations}
\begin{equation}
  d^{i}_{t+1} = \alpha_t^{i} d^{i}_{t} + \beta_t^{i} + \varepsilon_t^{i} \eqfinv
\end{equation}
where the non-stationary coefficients~$(\alpha^i_t, \beta^i_t)$ are, for all time~$t$,
solutions of the least-square problem
\begin{equation}
  (\alpha^i_t , \beta_t^i) = \argmin_{a,b}\sum_{s=1}^{N_{opt}}
  \norm{d^{i,s}_\post - ad^{i, s}_{t} - b}^2_2 \eqfinp
\end{equation}
The points $(d^{i, 1}_t, \ldots
d_t^{i, N_{opt}})$ correspond to the optimization scenarios.
The AR residuals $(\varepsilon_t^{el}, \varepsilon_t^{hw})$ are a white noise process.
\end{subequations}

\paragraph{Updating the forecast online}
once the AR model is calibrated, we use it to update the forecast during
assessment (see~\S\ref{MPC}). The update procedure is threefold:
\begin{itemize}
  \item[i)] we observe the demands~$w_t = (d^{el}_t, d^{hw}_t)$ between time
    $t-1$ and~$t$,
  \item[ii)] we update the forecast~$\overline w_{t+1}$ at time~$t+1$ with the AR model
    \begin{equation*}
      \overline w_\post =
      \bp{\overline d^{el}_\post, \overline d^{hw}_\post} =
      \Bp{\alpha_t^{el} d^{el}_{t} + \beta_t^{el}, \;
      \alpha_t^{hw} d^{hw}_{t} + \beta_t^{hw} }
      \eqfinv
    \end{equation*}
  \item[iii)] we set the forecast between time~$t+2$ and~$\final$ by
    using the mean values of the optimization scenarios:
    \begin{equation*}
      \overline w_\tau = \dfrac{1}{N^{opt}} \sum_{i=1}^{N^{opt}} w_\tau^{i}
      \quad \forall \tau = t+2, \cdots, \final \eqfinp
    \end{equation*}
\end{itemize}
Once the forecast $(\fw_\post, \ldots, \fw_\final)$ is available,
it is fed into the MPC algorithm that
solves Problem~\eqref{eq:mpcproblem}.

\subsubsection{SDDP procedure}
even if electrical and thermal demands are naturally correlated in
time~\cite{widen2010high}, the SDDP algorithm only relies upon marginal distributions.

\paragraph{Building offline probability distributions for SDDP}
rather than fitting an AR model like done with MPC,
we use the optimization scenarios to build marginal probability
distributions $ \mu_t^{of}$ that will feed the SDDP procedure
in~\eqref{eq:sddpforward}-\eqref{eq:sddpbackward}.

We cannot directly consider the discrete empirical marginal probability derived from
all $N_{opt}$~scenarios, because the support size would be too large for SDDP.
This is why we use optimal quantization to map the~$N_{opt}$ optimization
scenarios to~$S$ representative points.
We use a Loyd-Max quantization scheme~\cite{lloyd1982least}
to obtain a discrete probability distribution: at each time~$t$,
we use the~$N_{opt}$ optimization scenarios to build a partition
$\Xi = (\Xi_1, \cdots, \Xi_S)$, where~$\Xi$ is
the solution of the optimal quantization problem
\begin{equation}
\min_{\Xi} \quad
  \sum_{s=1}^S \bgp{ \sum_{w_t^i \in \Xi_s} \norm{w_t^i - \widetilde w_t^s}_2^2 }
\end{equation}
where~$\widetilde w_t^s = \frac{1}{card(\Xi_s)} \sum_{w_t^i \in \Xi_s} w_t^i$
is the so-called centroid of~$\Xi_s$.
Then, we set for all time $t = 0, \cdots, \final$ the discrete offline distributions
$ \mu_t^{of} =  \sum_{s=1}^S p_s \delta_{\widetilde w_t^s} \eqfinv$
where $\delta_{\widetilde w_t^s}$ is the
Dirac measure at point~$\widetilde w_t^s$ and $p_s = card(\Xi_s) / N^{opt}$
is the associated probability weight.
We have chosen~$S=20$ to have enough precision.

\paragraph{Computing value functions offline}
then, we use these probability distributions as an input to compute a set of value functions
with the procedure described in~\S\ref{sec:sddp}.

\paragraph{Using the value functions online}
once the value functions have been computed by SDDP, we are able to
compute online decisions with Equation~\eqref{eq:sddppolicydiscrete}\footnote{
  In practice, the quantization size of $\mu_t^{on}$ is bigger than
those of $\mu_t^{of}$, to have a greater accuracy online}.
SDDP, on the contrary of MPC, does not update the online probability distribution
$\mu_t^{on}$ during assessment to consider the information brought by the
previous observations.

\subsubsection{Heuristic procedure}
we choose to compare the MPC and SDDP algorithms with a basic
 decision rule. This heuristic is as follows: the battery is charged
whenever the solar production $\PV$ is available, and discharged to fulfill
the demand if there remains enough energy in the battery;
the tank is charged ($\va f^h_t > 0$) if the tank energy $\va h_t$ is lower than $\va h_0$,
the heater $\va f^t_t$ is switched on when the temperature is below the setpoint $\overline{\theta^i_t}$
and switched off whenever the temperature is above the setpoint plus a given margin.

\subsection{Results}

\subsubsection{Assessing on different meteorological conditions}
we assess the algorithms on three different days, with different meteorological
conditions (see Table~\ref{tab:meteoday}). Therefore,
we use three distinct sets of $N_{sim}$~assessment scenarios of demands, one for
each typical day.

\begin{table}[!ht]
 \centering
 \begin{tabular}{lccc}
   & Date & Temp. ($^\circ C$) & PV Production (kWh) \\
   \hline
   Winter Day & February, 19th & 3.3  & 8.4 \\
   Spring Day & April, 1st     & 10.1 & 14.8   \\
   Summer Day  & May, 31st      & 14.1 & 23.3
  \end{tabular}
  \caption{Different meteorological conditions}
  \label{tab:meteoday}
\end{table}
These three different days corresponds to different heating needs. During \emph{Winter
day}, the heating is maximal, whereas it is medium during \emph{Spring day} and
null during \emph{Summer day}. The production of the solar panel varies accordingly.

\subsubsection{Comparing the algorithms performances}
during assessment, we use MPC (see~\eqref{eq:mpcproblem}) and
SDDP (see~\eqref{eq:sddppolicy}) strategies to compute online decisions along
$N_{sim}$ assessment scenarios. Then, we compare the average electricity bill obtained
with these two strategies and with the heuristic.
The assessment results are given in Table~\ref{tab:comparison}:
means and standard deviation~$\sigma$ are computed by Monte Carlo with the $N_{sim}$ assessment
scenarios; the notation $\pm$ corresponds to the interval $\pm 1.96
\dfrac{\sigma}{\sqrt{N_{sim}}}$, which is a 95\% confidence interval.

\begin{table}[!ht]
  \centering
  \begin{tabular}{lccc}
    \toprule
                      & SDDP            & MPC             & Heuristic \\
    \midrule
    Offline time      & 50~s            & 0~s             & 0~s          \\
    Online time       & 1.5~ms          & 0.5~ms          & 0.005~ms          \\
    \midrule
    Electricity bill (\euro) &                 &                 &           \\
    \midrule
    Winter day        & 4.38~$\pm$ 0.02 & 4.59~$\pm$ 0.02 & 5.55~$\pm$ 0.02          \\
    Spring day        & 1.46~$\pm$ 0.01 & 1.45~$\pm$ 0.01 & 2.83~$\pm$ 0.01          \\
    Summer day        & 0.10~$\pm$ 0.01 & 0.18~$\pm$ 0.01 & 0.33~$\pm$ 0.02          \\
  \end{tabular}
  \caption{Comparison of MPC, SDDP and heuristic strategies}
  \label{tab:comparison}
\end{table}

We observe that MPC and SDDP exhibit close performance,
and make better than the heuristic.
If we consider mean electricity bills,
SDDP achieves better savings
than MPC during \emph{Summer day} and \emph{Winter day},
but SDDP and MPC display similar performances during \emph{Spring day}.

\begin{figure}[!ht]
  \centering
  \includegraphics[width=6cm]{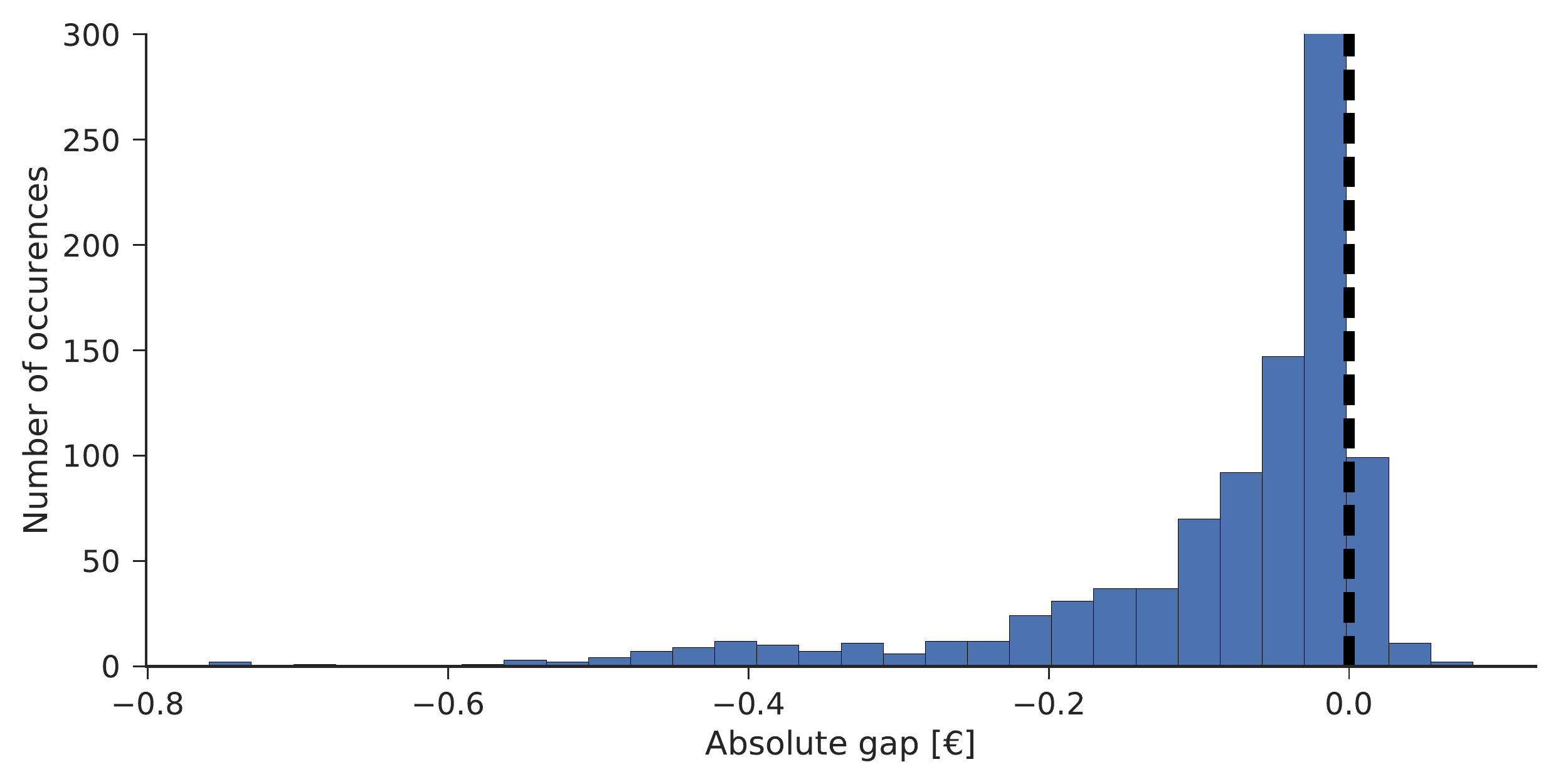}
  \caption{Absolute gap savings between MPC and SDDP during \emph{Summer day}}
  \label{fig:histsavings}
\end{figure}

In addition, SDDP achieves better savings than MPC for the vast majority of
scenarios. Indeed, if we compare the difference between
the electricity bills scenario by scenario, we observe that SDDP is better than
MPC for about 93\% of the scenarios. This can be seen on
Figure~\ref{fig:histsavings} that displays the histogram of the absolute gap savings
between SDDP and MPC during \emph{Summer day}.
The distribution of the gap exhibits a heavy tail that favors SDDP on
extreme scenarios.
Similar analyses hold for \emph{Winter} and \emph{Spring day}.
Thus, we claim that SDDP outperforms MPC for the electricity savings.
Concerning the performance on thermal comfort, temperature trajectories
are above the temperature setpoints specified in~\S\ref{sec:objective}
for both MPC and SDDP.

In term of numerical performance, it takes less than a minute to compute
a set of cuts as in~\S\ref{sec:sddp} with SDDP on a particular day. Then, the online computation
of a single decision takes 1.5~ms on average, compared to 0.5~ms for MPC.
Indeed, MPC is favored
by the linearity of the optimization Problem~\eqref{eq:mpcproblem},
whereas, for SDDP, the higher the quantization size~$S$, the slowest is the resolution
of Problem~\eqref{eq:sddppolicy}, but
the more information the online probability distribution~$\mu_t^{on}$ carries.

\subsubsection{Analyzing the trajectories}
we analyze now the trajectories of stocks in assessment, during \emph{Summer day}.
The heating is off, and the production of the solar panel is nominal at midday.

\begin{figure}[!ht]
  \centering
  \includegraphics[width=7cm]{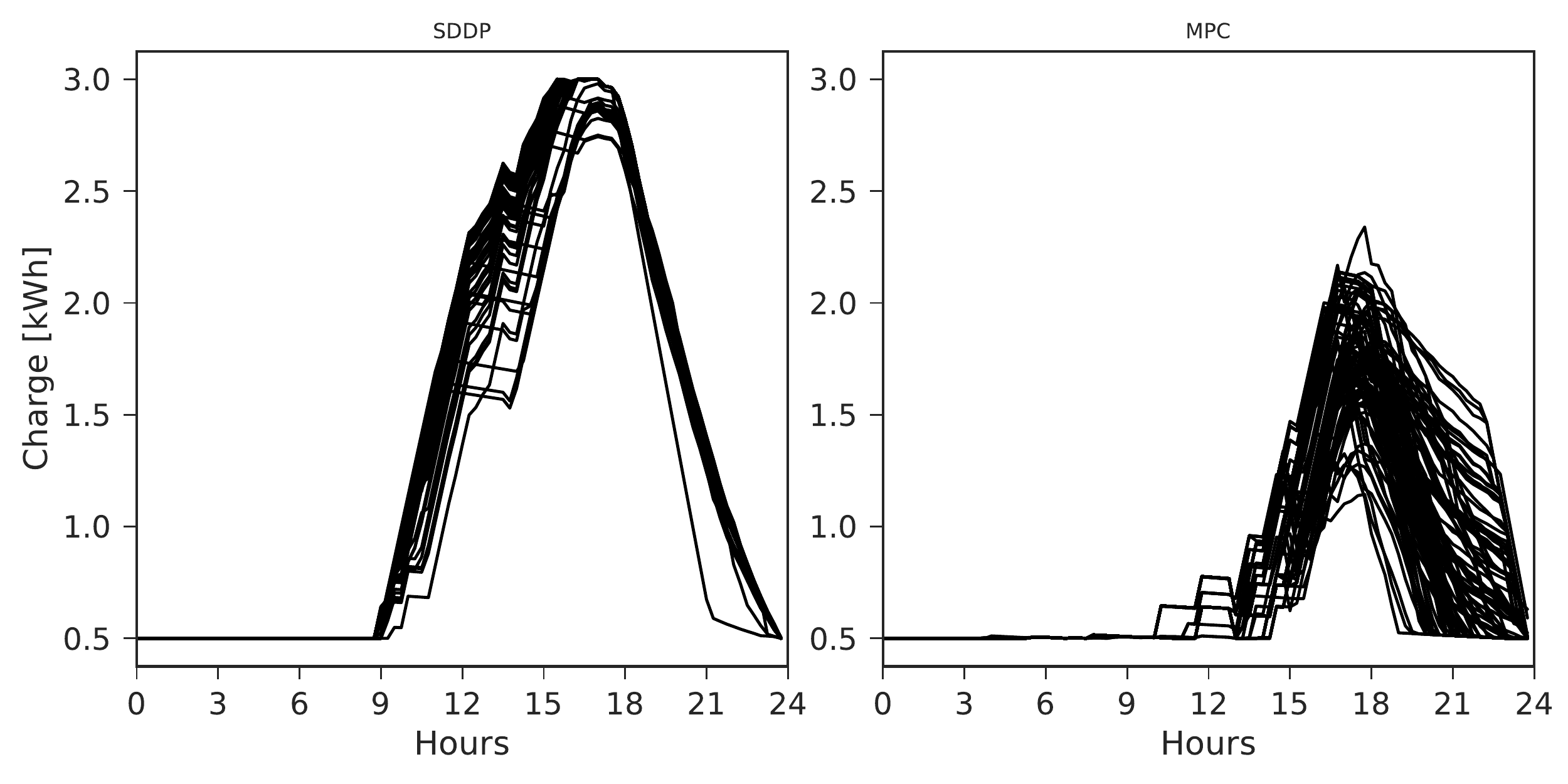}
  \caption{Battery charge trajectories for SDDP and MPC during \emph{Summer day}}
  \label{fig:battraj}
\end{figure}

Figure~\ref{fig:battraj} displays the state of charge of the battery along
a subset of assessment scenarios, for SDDP and MPC. We observe that SDDP
charges earlier the battery at its maximum.
On the contrary MPC charges
the battery later, and does not use the full potential of the battery. The
two algorithms discharge the battery to fulfill the evening demands.
We notice that each trajectory exhibits a single cycle of charge/discharge,
thus decreasing battery's aging.

\begin{figure}[!ht]
  \centering
  \includegraphics[width=7cm]{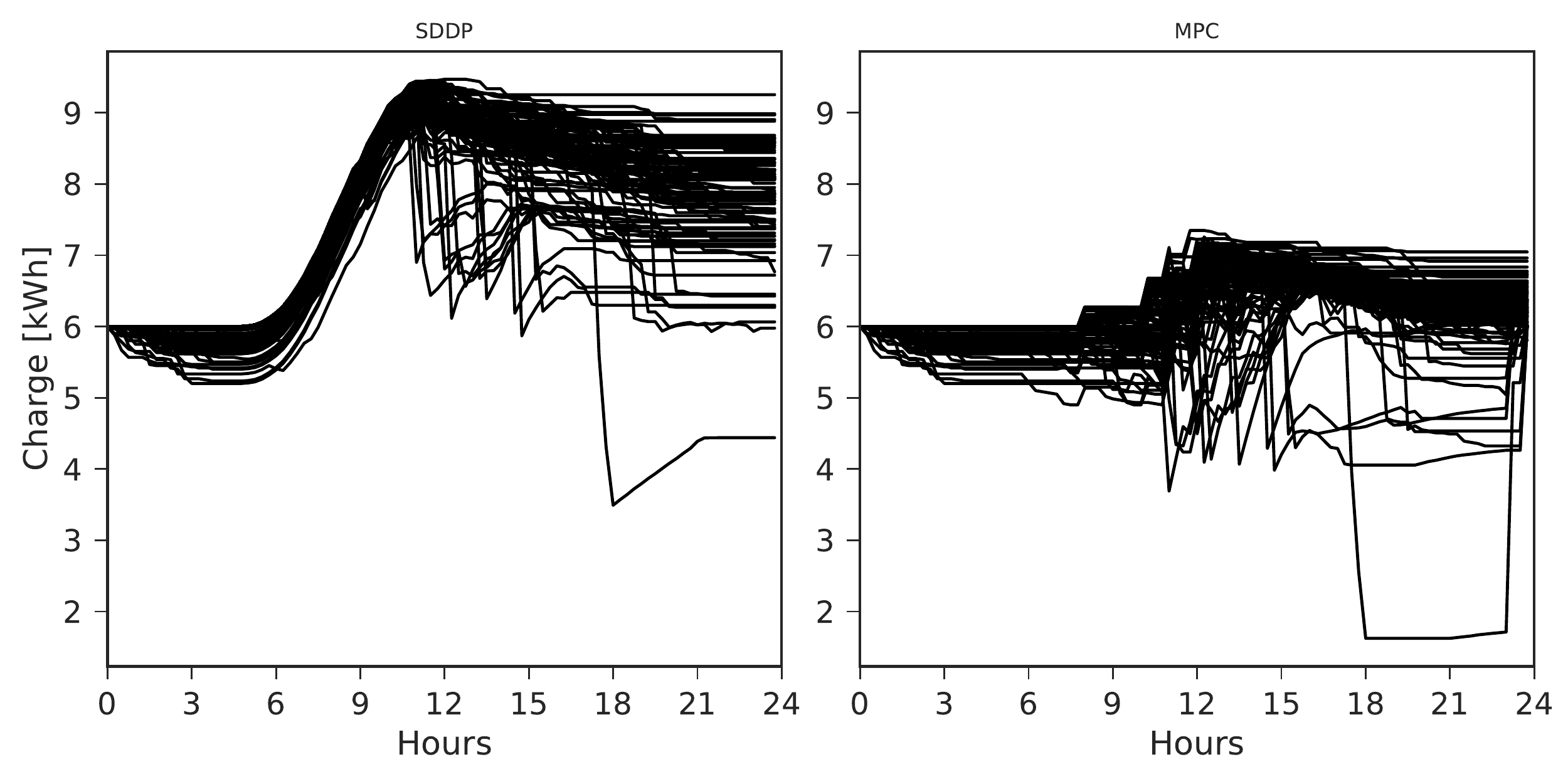}
  \caption{Hot water tank trajectories for SDDP and MPC during \emph{Summer day}}
  \label{fig:dhwtraj}
\end{figure}

Figure~\ref{fig:dhwtraj} displays the charge of the domestic hot water tank along
the same subset of assessment scenarios. We observe
a similar behavior as for the battery trajectories:
SDDP uses more the electrical hot water tank to store the excess of PV energy,
and the level of the tank is greater at the end of the day than in MPC.

This analysis suggests that SDDP makes a better use of storage capacities than
MPC.

\section{Conclusion}

We have presented a domestic microgrid energy system,
and compared different optimization algorithms to control
the stocks with an Energy Management System.

The results show that optimization based strategies outperform the proposed
heuristic procedure in term of money savings.
Furthermore, SDDP outperforms MPC during \emph{Winter} and \emph{Summer day} ---
achieving up to 35\% costs savings --- and displays similar performance as MPC
during \emph{Spring day}. Even if SDDP and MPC exhibit close performance,
a comparison scenario by scenario shows that SDDP beats MPC
most of the time (more than 90\% of scenarios during \emph{Summer day}).
Thus, we claim that SDDP is better than MPC to manage uncertainties
in such a microgrid, although MPC gives also good performance.
SDDP also makes a better use of storage capacities.

Our study can be extended in different directions.
First, we could mix SDDP and MPC to recover the benefits of these two algorithms.
Indeed, SDDP is designed to handles the uncertainties' variability but fails to capture
the time correlation, whereas MPC ignores the uncertainties'
variability, but considers time correlation by means of
a multistage optimization problem.
Second, we are currently investigating the optimization of
larger scale microgrids --- with different interconnected buildings ---
by decomposition methods.

\end{document}